\documentclass[12pt]{article}
\usepackage{amsmath,amssymb,amsthm,amscd}
\usepackage{diagram}
\title{A survey on  symplectic singularities and resolutions}
\author{Baohua Fu}
\newtheorem{Thm}{Theorem}[section]

\newtheorem{Prop}[Thm]{Proposition}
\newtheorem{Cor}[Thm]{Corollary}
\newtheorem{Conj}{Conjecture}
\newtheorem{Rque}{Remark}[section]

\newtheorem{Exam}[Thm]{Example}
\newtheorem{Pb}{Problem}

\def\cit{{\mathbb C}}

\def\qit{{\mathbb Q}}
\def\zit{{\mathbb Z}}
\def\pit{{\mathbb P}}
\def\0{{\mathcal O}}
\def\g{{\mathfrak g}}
\def\h{{\mathfrak h}}
\def\p{{\mathfrak p}}

\def\n{{\mathfrak n}}
\def\X{{\mathcal X}}
\def\Y{{\mathcal Y}}
\def\u{{\mathfrak u}}

\begin{document}
\maketitle
\tableofcontents

\section*{Introduction}

This is a survey written in an expositional style  on the topic of
symplectic singularities and symplectic resolutions, which could
also serve as an introduction to this subject.

We work over the complex number field.
 A normal variety $W$ is called a {\em symplectic
variety} if its smooth part admits a holomorphic symplectic form
$\omega$ whose pull-back to any resolution $\pi: Z \to W$ extends
to a holomorphic 2-form $\Omega$ on $Z$. If furthermore the
extended 2-form $\Omega$ is a symplectic form, then $\pi$ is
called a {\em symplectic resolution}.

The existence and non-existence of symplectic resolutions are
difficult to decide. However one hopes that a symplectic variety
admits at most finitely many non-isomorphic symplectic resolutions
(section \ref{finite}).

Symplectic resolutions behave much like hyperK\"ahler manifolds.
Motivated by the work of D. Huybrechts (\cite{Huy}), one expects
that two symplectic resolutions are deformation equivalent (section
\ref{deformation}). This would imply the invariance of the
cohomology of the resolution, which is expected to be recovered by
the Poisson cohomology of the symplectic variety (section
\ref{cohomology}). As a special case of Bondal-Orlov's conjecture,
one expects that two symplectic resolutions are derived equivalent
(section \ref{derived}).  Finally motivated by the results in
dimension 4, one expects some simple birational geometry in
codimension 2 for symplectic resolutions (section \ref{birational}).

Examples of symplectic varieties include  quotients $\cit^{2n}/G$
with $G$ a finite subgroup of $Sp(2n)$ and normalizations of
nilpotent orbit closures in semi-simple Lie algebras.

For symplectic resolutions of nilpotent orbit closures, our
understanding is more or less complete. However, our knowledge of
symplectic resolutions of quotient singularities is rather limited
(only cohomology and derived equivalence have been fully
understood).

A funny observation is that all known examples of symplectic
resolutions are modeled locally on Hilbert schemes or on Springer's
resolutions. For O'Grady's symplectic resolution of moduli space of
sheaves (\cite{O'G1}, \cite{O'G2}), Kaledin and Lehn (\cite{KL}) proved that it is modeled
locally on a Springer's resolution of a nilpotent orbit closure in
$\mathfrak{sp}(4)$ (see also \cite{LS}). It would be very interesting to find out more
local models. \vspace{0.4cm}

{\em Acknowledgement:} This paper is by no means intended to be a
complete account of the topic. Therefore I apologize for the
omission of any relevant work and references. I would like to take
this opportunity to thank J. Alev and T. Lambre for their
invitation to the journey ``Alg\`ebre et G\'eom\'etrie de Poisson"
(Clermond-Ferrand 2004), where they encouraged me to write this
survey.  I would like to thank D. Kaledin and  J. Sawon for
corrections and suggestions to a preliminary version of this
paper.

\section{Basic definitions and properties}
\subsection{symplectic singularities}
Since A. Beauville's pioneering paper \cite{Be3}, symplectic
singularities have received a particular attention by many
mathematicians. As explained in \cite{Be3}, the motivation of this
notion comes from the analogy between rational Gorenstein
singularities and Calabi-Yau manifolds.

By a {\em resolution} we mean a proper surjective morphism $\pi: Z
\to W$ such that: (i). $Z$ is smooth; (ii). $\pi^{-1}(W_{reg}) \to
W_{reg}$ is an isomorphism. If furthermore $\pi$ is a projective
morphism, then $\pi$ is called a {\em projective resolution}.

Recall that a compact K\"ahler manifold of dimension $m$ is a
Calabi-Yau manifold if it admits a nowhere vanishing holomorphic
$m$-form. Its singular counterpart is varieties with rational
Gorenstein singularities, i. e.  normal varieties $W$ of dimension
$m$ whose smooth part  admits  a holomorphic nowhere vanishing
$m$-form such that its pull-back to any resolution $Z \to W$
extends to a holomorphic form on the whole of $Z$.

A holomorphic 2-form on a smooth variety is called {\em
symplectic} if it is closed and non-degenerate at every point.
Among Calabi-Yau manifolds, there are symplectic manifolds, i.e.
smooth varieties which admit a holomorphic symplectic form. By
analogy, the singular counterpart of symplectic manifolds is
varieties with symplectic singularities (also called {\em
symplectic varieties}), i.e. normal varieties $W$ whose smooth
part admits  a holomorphic symplectic form $\omega$ such that its
pull-back to any resolution $Z \to W$ extends to a holomorphic
2-form $\Omega$ on $Z$.

One should bear in mind that the 2-form $\Omega$ is closed but may
be degenerated at some points. The following proposition follows
from a theorem of Flenner on extendability of differential forms.
\begin{Prop}
A normal variety with singular part having codimension $\geq 4$ is
symplectic if and only if its smooth part admits a holomorphic
symplectic form.
\end{Prop}

More generally, one has the following characterization of symplectic
varieties:
\begin{Thm}[Namikawa \cite{Nam2}] \label{ratGor}
A normal variety is symplectic if and only if it has only rational
Gorenstein singularities and its smooth part admits a holomorphic
symplectic form.
\end{Thm}

\subsection{stratification theorem}
In differential geometry, it is well-known that every Poisson
structure on a real smooth manifolds gives rise to a foliation by
symplectic leaves. The following stratification theorem extends this
to the case of symplectic varieties.
\begin{Thm}[Kaledin \cite{Ka4}]
Let $W$ be a symplectic variety. Then there exists a canonical
stratification $W= W_0 \supset W_1 \supset W_2 \cdots $ such that:

(i). $W_{i+1}$ is the singular part of $W_i$;

(ii). the normalization of every irreducible component of $W_i$ is a
symplectic variety.
\end{Thm}

One shows easily that every $W_i$ is a Poisson subvariety  in $W$.
The difficult part is to show that the normalization of any Poisson
subscheme in $W$ is still a symplectic variety (Theorem 2.5
\cite{Ka4}). An immediate corollary is

\begin{Cor}\label{evencodim}
Every irreducible component of the singular part of a symplectic
variety has even codimension.
\end{Cor}
It has been previously proved  by Y. Namikawa (\cite{Nam3}) that the
singular part of a symplectic variety has no codimension 3
irreducible components.

\begin{Cor}
Let $W$ be a symplectic variety  which is locally a complete
intersection, then the singular locus of $W$ is either empty or  of
pure codimension 2.
\end{Cor}

In fact it is proved in \cite{Be3} (Proposition 1.4) that
$Sing(W)$ is of codimension $\leq 3$. Now Corollary
\ref{evencodim} excludes the case of codimension 3. This corollary
gives rise to the following
\begin{Pb}
 Classify symplectic varieties which are of complete intersection, and those which admit a
symplectic resolution.
\end{Pb}

Such examples include nilpotent cones in semi-simple Lie algebras
and rational double points ($ADE$ singularities).

\subsection{symplectic resolutions}

By Hironaka's big  theorem, any complex variety admits a
resolution, but there may exist many different resolutions.   One
would like to find out some ``good" resolutions. In dimension 1,
the resolution is unique, which is given by the normalization. In
dimension 2, one also finds  a ``good" resolution, the so-called
minimal resolution, i.e. any other resolution factorizes through
this resolution. When the dimension is higher, one finds the
following class of  preferred resolutions (crepant resolutions).
However this notion is defined only for varieties with a
$\qit$-Cartier canonical divisor.

Let $W$ be a normal variety. A Weil divisor $D$ on $W$ is called
{\em $\qit$-Cartier} if some non-zero multiple of $D$ is a Cartier
divisor.  $W$ is called {\em $\qit$-factorial} if every Weil divisor
on $W$ is {\em $\qit$-Cartier}. The quotient of a smooth variety by
a finite group is $\qit$-factorial.

For a normal variety $W$, the closure of a canonical divisor of
$W_{reg}$ in $W$ is called a {\em canonical divisor} of $W$, denoted
by $K_W$. In general it is only a Weil divisor. Suppose that $K_W$
is $\qit$-Cartier, i.e. there exists some non-zero integer $n$ such
that $\0_W(n K_W)$ is an invertible sheaf. Then for any resolution
$\pi: Z \to W$, the pull-back $\pi^*(K_W):= \frac{1}{n} \pi^*
(nK_W)$ is well-defined. The resolution $\pi$ is called {\em
crepant} if $K_Z \equiv \pi^*(K_W)$, i.e. if $\pi$ preserves the
canonical divisor.

One should bear in mind that a crepant resolution does not exist
always (as we will show soon). However whenever such a  resolution
exists, there is a close relationship between the geometry of the
resolution and the geometry of the singular  variety.

For a symplectic variety $W$, a resolution $\pi: Z \to W$ is
called {\em symplectic} if the extended 2-form $\Omega$ on $Z$ is
non-degenerate, i.e. if $\Omega$ defines a symplectic structure on
$Z$. One might wonder if this definition depends on the special
choice of the symplectic structure on the smooth part. However we
have
\begin{Prop}
Let $W$ be a symplectic variety and $\pi: Z \to W$ a resolution.
Then the following statements are equivalent:

(i). $\pi$ is crepant;

(ii). $\pi$ is symplectic;

(iii). the canonical divisor $K_Z$ is trivial;

(iv). for any symplectic form $\omega'$ on $W_{reg}$, the pull-back
$\pi^* (\omega')$ extends to a symplectic form on $Z$.
\end{Prop}
\begin{proof}
The only implication to be proved is $(i) \Rightarrow (iv)$: Since
$W$ is symplectic, it has only rational Gorenstein singularities.
Now by \cite{Nam2}, any symplectic form $\omega'$ on $W_{reg}$
extends to a holomorphic 2-form $\Omega'$ on $Z$. Notice that
$\wedge^{top} \omega'$ gives a trivialization of $K_W$,
$\wedge^{top} \Omega'$ also trivializes $K_Z$, since $\pi$ is
crepant, thus $\Omega'$ is symplectic.
\end{proof}

\subsection{Namikawa's work}

 We now give a quick review of Y. Namikawa's important work
on symplectic varieties, while some of his work on symplectic resolutions
will be explained in detail later.

  In \cite{Nam2} (Theorem 7 and 8),
he  proved a stability theorem and a local Torelli theorem. The
deformation theory of  such a variety $W$ is studied in
\cite{Nam1} (Theorem 2.5), where he proved that if $codim\ Sing(W)
\geq 4$, then the Kuranishi space $Def(W)$ is smooth.

In \cite{Nam3}, it is shown that a 
symplectic variety has terminal
singularities if and only if the codimension 
of the singular part is at least 4. Any flat deformation of 
a projective $\qit$-factorial terminal symplectic variety 
is locally trivial, i.e. it is not smoothable by flat deformations
(see \cite{Nam6}). Notice that such a variety admis no 
symplectic resolutions (see the proof of Proposition \ref{purecodim2}).

It is conjectured in \cite{Nam6} that a projective symplectic variety 
is smoothable by flat deformations if and only if it admits 
a symplectic resolution.   

Recently a local version of \cite{Nam6} is obtained by Y.  Namikawa
himself in \cite{Nam7}. A symplectic variety is called {\em convex} if
there exists a projective birational morphism to an affine normal variety.
For convex symplectic varieties, there exists a theory
 of Poisson deformations (see for example \cite{GK}). The main theorem of
\cite{Nam7} says that for a convex symplectic variety with terminal 
singularities $W$ such that $W^{an}$ is $\qit$-factorial, 
any Poisson deformation of $X$ is locally trivial (forgetting the Poisson
structure).

\section{Examples}
\subsection{quotient singularities}
Let $W$ be a quasi-projective symplectic variety and $G$ a finite
subgroup of $Aut(W)$ preserving a symplectic form on $W_{reg}$.
The symplectic form on $W^0:=W_{reg} - \cup_{g \neq 1} Fix(g)$
descends to a symplectic form on $W^0/G$, which extends to a
symplectic form on $(W/G)_{reg}$, since the complement of $W^0/G$
in $W/G$ has codimension $\geq 2$.  Now it is shown in \cite{Be3}
(Proposition 2.4) that this symplectic form extends to a
holomorphic 2-form in any resolution. In conclusion the quotient
$W/G$ is a symplectic variety. Here are some special cases:
\begin{Exam}
Let $G$ be a finite sub-group of $SL(2, \cit)$. The quotient
$\cit^2/G$ is a symplectic variety with rational double points. It
admits a unique symplectic resolution, given by the minimal
resolution.
\end{Exam}
\begin{Exam}
Let $V$ be a finite-dimensional symplectic vector space and $G <
Sp(V)$ a finite sub-group. Then the quotient $V/G$ is a symplectic
variety. However, it is difficult to decide when $V/G$ admits a
symplectic resolution.
\end{Exam}
\begin{Exam}
Let $Y$ be a smooth quasi-projective variety and $G < Aut(Y)$ a
finite group. Then $G$ acts on $T^*Y$ preserving the natural
symplectic structure, thus the quotient $(T^*Y)/G$ is a symplectic
variety.
\end{Exam}
\begin{Exam}
Let $W$ be a symplectic variety. Then the symmetric product
$W^{(n)}$ is a symplectic variety.  When $W$ is smooth and $dim(W)
\geq 4$,  $W^{(n)}$ does not admit any symplectic resolution
(see Proposition \ref{purecodim2}).  However when $W$ is a smooth
symplectic surface $S$ (for example an Abelian surface, a K3 surface
or the cotangent bundle of a curve),  $S^{(n)}$ admits a
symplectic resolution given by $Hilb^n(S) \to S^{(n)}$.  This is
also the unique projective symplectic resolution  of $S^{(n)}$ (see
\cite{FN}).
\end{Exam}
\begin{Exam} \label{mukai}
Let $G$ be a finite subgroup of $SL(2)$ and $S \to \cit^2/G$  the
minimal resolution. Then the symmetric product $(\cit^2/G)^{(n)}$ is
naturally identified with $\cit^{2n}/G'$, where $G'$ is the wreath
product of $S_n$ with $G$. Now a symplectic resolution of
$\cit^{2n}/G'$ is given by the composition $$ Hilb^n(S) \to S^{(n)}
\to Sym^n(\cit^2/G).$$

Consider the case $G = \pm 1 $. Then $S = T^* \pit^1$. The central
fiber of the symplectic resolution $\pi: Hilb^n(S) \to
(\cit^2/G)^{(n)}$ contains a component isomorphic to $\pit^n$. By
performing a Mukai flop (for details see section \ref{birational})
along this component, one obtains another symplectic resolution
which is non-isomorphic to $\pi$. More discussions on this example
can be found in \cite{Fu5}.
\end{Exam}

\subsection{nilpotent orbit closures}

Let $\g$ be a  semi-simple complex Lie algebra,
 i.e. the bilinear form $\kappa (u, v) := trace (ad_u
\circ ad_v)$ is non-degenerate, where $ad_u: \g \to \g$ is the
linear map given by $x \mapsto [u,x]$. Let $Aut(\g) = \{\phi \in
GL(\g) | [\phi(u), \phi(v)] = [u,v], \forall u, v \in \g \}$,
which is a Lie group but may be disconnected, whose identity
connected component is the adjoint group $G$ of $\g$.

 An element $v \in \g$ is called {\em
semi-simple} (resp. {\em nilpotent}) if the linear map $ad_v$ is
semi-simple (resp. nilpotent), whose orbit under the natural action
of $G$ is denoted by $\0_v$, which is called a semi-simple orbit
(resp. nilpotent orbit).

Semi-simple orbits in $\g$ are parameterized by $\h / W$, where
$\h$ is a Cartan sub-algebra in $\g$ and $W$ is the Weyl group. In
particular there are infinitely many semi-simple orbits in $\g$.
Semi-simple orbits possed a rather simple geometry, for example
they are closed and simply-connected.

To the contrary, nilpotent orbits have a much richer geometry. The
classification of nilpotent orbits has been carried out by
Kostant, Dynkin, Bala-Carter et. al. via either weighted Dynkin
diagrams or partitions in the case of classical simple Lie
algebras.

\begin{Exam}
An element in $\mathfrak{sl}_{n+1}$ is nilpotent if and only if it
is conjugate to some matrix $diag(J_{d_1}, \ldots, J_{d_k})$, where
$J_{d_i}$ is a $d_i \times d_i$ Jordan bloc with zeros on the
diagonal, and $d_1 \geq d_2 \geq \cdots \geq d_k \geq 1$ are
integers such that $\sum_{i=1}^k d_i = n+1$, i.e. $[d_1, \ldots,
d_k]$ is a partition of $n+1$. This gives a one-one correspondence
between nilpotent orbits in $\mathfrak{sl}_{n+1}$ and partitions of
$n+1$.
\end{Exam}

%\begin{Thm} [\cite{CM}]
%(i). Nilpotent orbits in $\mathfrak{so}_{2n+1}$ are bijective to
%partitions of $2n+1$ such that even parts occur with even
%multiplicity;

%(ii). Nilpotent orbits in $\mathfrak{sp}_{2n}$ are bijective to
%partitions of $2n$ such that odd parts occur with even
%multiplicity;

%(iii). Nilpotent orbits in $\mathfrak{so}_{2n}$ are bijective to
%partitions of $2n$ such that even parts occur with even
%multiplicity, except ``very even" partitions, i.e. those with only
%even parts and each with even multiplicity. Every very even
%partition corresponds to two nilpotent orbits.
%\end{Thm}

For other classical simple Lie algebras, a similar description of
nilpotent orbits exists (see \cite{CM}).  The following theorem is
fundamental in the study of nilpotent orbits.

\begin{Thm}[Jacobson-Morozov] \label{JM}
Let $\g$ be a complex semi-simple Lie algebra and $v \in \g$ a
nilpotent element. Then there exist two elements $H, u \in \g$
such that $[H, v] = 2v,  [H, u] = -2u, [v, u] = H.$
\end{Thm}

The triple $\{H, v, u \}$ is called a {\em standard triple}, which
provides an isomorphism $\phi: \mathfrak{sl}_2 \to  \cit \langle H,
v, u \rangle $. Now $\g$ becomes an $\mathfrak{sl}_2$-module via
$\phi$ and the adjoint action. Thus $\g$ is decomposed as $\g =
\oplus_{i \in \zit} \g_i$, where $\g_i = \{x \in \g \  | \ [H, x] =
i x \}.$ Let $\p = \oplus_{i \geq 0} \g_i$ and $P$ a connected
subgroup of $G$ with Lie algebra $\p$. Then $\p$ is a parabolic
sub-algebra of $\g$ and $P$ is a parabolic subgroup of $G$. Let $\n
= \oplus_{i \geq 2} \g_i$  and $\u = \oplus_{i \geq 1} g_i$. The
nilpotent orbit $\0_v$ is called {\em even} if $\g_{1} = 0$ or
equivalently if $\g_{2k+1} = 0$ for all $k \in \zit$. In this case,
one has $\u = \n \simeq (\g/\p)^*$.

The nilpotent orbit $\0_v$ is not closed in $\g$. Its closure
$\overline{\0_v}$ is a singular (in general non-normal) variety.
 There exists a natural
proper resolution of $\overline{\0_v}$ as follows (the so-called
{\em Springer's resolution}): $G \times^P \n \xrightarrow{\pi}
\overline{\0_v},$ where $ G \times^P \n$ is the quotient of $G
\times \n$ by $P$ acting as $p(g, n) = (gp^{-1}, Ad_p(n))$. The
group $G$ acts on $ G \times^P \n$ by $g (h, n) = (gh, n)$. Then the
resolution is $G$-equivariant and maps the orbit $G\cdot (1, v)$
isomorphically to $\0_v$.

For any $g \in G$, the tangent space $T_{Ad(g)v} \0_v$ is
isomorphic to $[\g, Ad(g)v]$. Now we define a 2-form $\omega$ on
$\0_v$ as follows: $$ \omega_{Ad(g)v} ([u_1, Ad(g)v], [u_2, Ad(g)
v]) = \kappa([u_1, u_2], Ad(g)v).$$  The 2-form $\omega$ is in
fact a closed non-degenerate 2-form, i.e. a holomorphic symplectic
form on $\0_v$ (the so-called Kostant-Kirillov-Souriau form).

\begin{Prop}(\cite{Be2}, \cite{Pan})
The symplectic form $\pi^* \omega$ on $G\cdot (1, v)$ can be
extended to a global 2-form $\Omega$ on $G \times^P \n$.
\end{Prop}
\begin{proof}
Take an element $(g, n) \in G \times^P \n$, then the tangent space
of $G \times^P \n$ at $(g, n)$ is canonically isomorphic to $\g
\times \n / \{(x, [n,x]) | x \in \p \}.$ We define a 2-form
$\beta$ on $G \times^P \n$ as follows:
$$ \beta_{(g,n)}((u, m), (u', m')) = \kappa([u, u'],n)+\kappa(m',u) -
\kappa(m, u').$$ The kernel of $\beta_{(g,n)}$ is $\{(u, [n,u]) \ |
\ u \in \oplus_{i \geq -1} \g_i \}$, which contains the sub-space
$\{(x, [n,x]) | x \in \p \}$, thus this 2-form descends to a 2-form
$\Omega$ on $G \times^P \n$. The 2-forms $\Omega$ and $\pi^* \omega$
coincide at the point $(1, v)$, thus they coincide on $G(1, v)$
since both are $G$-equivariant.
\end{proof}
\begin{Cor}
The 2-form $\Omega$ is symplectic if and only if $\g_{-1} = 0$,
i.e. if and only if $\0_v$ is an even nilpotent orbit.
\end{Cor}
\begin{proof}
The kernel of $\Omega$ at $(g, 0)$ is isomorphic to $\g_{-1}$,
thus if $\Omega$ is symplectic, then $\g_{-1} = 0$. Conversely if
$\g_{-1} = 0$, then $G \times^P \n = G \times^P \u \simeq T^*
(G/P)$, which implies that the canonical bundle  $K$ of $G
\times^P \n$ is trivial. Notice that $\Omega^{top}$ gives a
trivialization  of $K$, thus $\Omega^{top}$ is non-zero
everywhere, i.e. $\Omega$ is symplectic.
\end{proof}

Notice that the resolution $G \times^P \n \xrightarrow{\pi}
\overline{\0}$ factorizes through the normalization $\widetilde{\0}
\to \overline{\0}$, which gives a resolution of $\widetilde{\0}$.
\begin{Cor}\label{nilpsym}
The normalization $\widetilde{\0}$ of a nilpotent orbit closure in a
complex semi-simple Lie algebra is a symplectic variety. The
resolution $G \times^P \n \to \widetilde{\0}$ is symplectic if and
only if $\0$ is an even nilpotent orbit.
\end{Cor}

One should remember  that even for an even nilpotent orbit
closure, there can exist some symplectic resolutions not of the
above form.

\subsection{isolated singularities}

Let $V$ be a finite-dimensional symplectic vector space and $G <
Sp(V)$ a finite sub-group. Suppose furthermore that the non trivial
elements in $G$ have all their eigenvalues different to 1, then the
quotient $G/V$ is a symplectic variety with an isolated singularity,
which admits a symplectic resolution if and only if $dim(V) =2$ (see
Corollary \ref{isosing}).

For example, let $\xi$ be the primitive cubic unit root, which acts
on $\cit^{2n}$ by the multiplication of $\xi $ on the first $n$
coordinates and by the multiplication of $\xi^2 $ on the last $n$
coordinates. Then the quotient $\cit^{2n}/\langle \xi \rangle$ is an
isolated symplectic singularity. A characterization of this
singularity has been given in \cite{Dru}.

Another type of isolated symplectic singularities comes from minimal
nilpotent orbits $\0_{min}$, i.e. the unique non-zero nilpotent
orbit which is contained in the closure of all non-zero nilpotent
orbits. Then $\overline{\0}_{min} = \0_{min} \cup \{0\}$ is normal
with an isolated symplectic singularity.  Conversely an isolated
symplectic singularity with smooth projective tangent cone is
analytically isomorphic to $\overline{\0}_{min}$ (see \cite{Be3}).
It is suggested in \cite{Be3} to classify isolated symplectic
singularities with trivial local fundamental group.

Among minimal nilpotent orbit closures, only those in $sl(n+1)$
admit a symplectic resolution (see Proposition \ref{omin}). In
this case, $\0_{min}$ consists of matrices of trace zero and rank
1. A symplectic resolution is given by $T^* \pit^n \to
\overline{\0}_{min}$. It is believed  in  that a projective
isolated symplectic singularity (of dimension $\geq 4$) admitting
a symplectic resolution is isomorphic to a such singularity (see
Corollary \ref{isosing} for the case of quotient isolated
singularities). Some discussions are given in \cite{CMSB}.

Here is a deformation of this symplectic resolution. In the
following, a point in $\pit^n$ will also be thought of a line in
$\cit^{n+1}$. Let $$ \mathcal{Z} = \{(l, A, a) \in \pit^n \times
\mathfrak{gl}_{n+1} \times \cit | Im(A) \subset l; Av = av, \
\forall v \in l\}$$ and $\mathcal{W}' = \{(A, a) \in
\mathfrak{gl}_{n+1} \times \cit | A^2 = a A; rk(A) = 1 \} $. We
denote by $\mathcal{W}$ the closure of $\mathcal{W}'$ in
$\mathfrak{gl}_{n+1} \times \cit$. Then the natural map $\mathcal{Z}
\xrightarrow{f} \mathcal{W}$ is a deformation of the symplectic
resolution $T^* \pit^n \to \overline{\0}_{min}$. Notice that for $a
\neq 0$ the map between fibers $\mathcal{Z}_a \xrightarrow{f_a}
\mathcal{W}_a$ is an isomorphism.

\section{Semi-smallness}

Recall that a morphism $\pi: Z \to W$ is called {\em semi-small}
if for every closed subvariety $F$ in $Z$, we have $2 \ codim\  F
\geq codim \ (\pi(F))$. This is a remarkable property, which
enables us, for example, to use the intersection cohomology
theory.

\begin{Exam}
Let $S$ be a normal surface, then any resolution of $S$ is
semi-small. However this is not the case in higher dimension. In
fact, the blowup of a point in the exceptional locus of a resolution
gives a resolution which is never semi-small.
\end{Exam}

As discovered partially by J. Wierzba (\cite{Wie}), Y. Namikawa
(\cite{Nam2}), and Hu-Yau (\cite{HY}), then in full generality by
D. Kaledin (\cite{Ka3}, \cite{Ka4}), symplectic resolutions enjoy
the semi-small property.
\begin{Thm}\label{semismall}
Let $W$ be a normal algebraic variety with only rational
singularities and $\pi: Z \to W$ a proper resolution. Suppose that
$Z$ admits a symplectic structure $\Omega$. Then the resolution
$\pi$ is semi-small.  In particular, a  symplectic resolution is
semi-small.
\end{Thm}
\begin{Rque}
When $\pi$ is projective, this theorem has been proved by D.
Kaledin (\cite{Ka4}, Lemma 2.7). With minor changes, his proof
works also for the proper case.
\end{Rque}
\begin{proof}
Let $Y \subset W$ be an irreducible closed subvariety and $F$ an
irreducible component of $\pi^{-1}(Y)$. One needs to show that $2
codim\ (F) \geq codim \  \pi(F)$.

By Chow's lemma (see for example \cite{Har}, Chap. II, exercise
4.10), there exists an algebraic variety $F'$ and a birational
proper morphism $f: F' \to F$ such that the composition $\pi \circ
f: F' \to Y$ is a projective morphism.

Now take a projective resolution $X \to F'$ and denote by $\sigma$
the composition morphism from $X$ to $Y$. Let $\eta: X \to Z$ be the
composition. By shrinking $W$ and $Y$ if necessary, we can assume
that

(i). $W$ is affine;

(ii). $Y$ is smooth and $Y = \pi(F) = \sigma (X)$;

(iii). $\sigma$ is smooth (this is possible since $X$ is smooth, see
\cite{Har} Corollary 10.7,Chap. II).
$$
$$

 For any $y \in Y$, we denote by $X_y$ (resp. $F'_y, F_y$)
the fiber over $y$ of the morphism $\sigma$ (resp. $\pi \circ f,
\pi$).

The arguments in \cite{Ka4} show that for any $x\in X_y$, we have
$T_x X_y \subset ker (\eta^* \Omega)_x$. Let $U \subset F$ be the
open set such that $\eta^{-1}(U) \to  U$ is an isomorphism. Then
for any point $z \in U$, $T_z F_y$ and $T_z F$ are orthogonal with
respect to $\Omega$. By our assumption, $\Omega$ is non-degenerate
everywhere on $Z$, thus
$$ dim (T_z F_y) + dim (T_z F) \leq dim (Z),$$ which gives the inequality in the theorem.
\end{proof}

\begin{Prop}\label{purecodim2}
Let $X$ be a smooth irreducible symplectic variety and $G$ a finite
group of symplectic automorphisms on $X$. Suppose that $V/G$ admits
a  symplectic resolution, then the subvariety $F= \cup_{g \neq 1}
Fix(g)\subset X $ is either empty or of pure codimension 2 in $X$.
\end{Prop}
\begin{proof}
Being a quotient of a $\qit$-factorial normal variety  by a finite
group, $V/G$ is again $\qit$-factorial and normal. This gives that
any component $E$ of the exceptional locus of a proper resolution
$\pi: Z \to X/G$ is of  codimension 1. On the other hand, if $\pi$
is a symplectic resolution, then by the semi-smallness, we have $
2 = 2 codim (E) \geq codim (\pi(E)).$ Suppose that $X/G$ is not
smooth, then the singular locus of $X/G$ is of codimension $\geq
2$, hence $codim(\pi(E)) = 2$.

However, the singular locus of $X/G$ is contained in $p(F)$, where
$p: X \to X/G$ is the natural map, hence $codim(F) \geq 2$. Notice
that for any $g\in G$, $Fix(g)$ is of even dimension since $g$ is
symplectic, thus $F$ has no codimension 1 component, i.e. $F$ is
of pure codimension 2.
\end{proof}

\begin{Cor}\label{isosing}
Let $X$ be a smooth irreducible symplectic variety and $G$ a finite
group of symplectic automorphisms on $X$. Suppose that $X/G$ has
only isolated singularities, then $X/G$ admits a  symplectic
resolution if and only if $dim(X) = 2$.
\end{Cor}

\section{Quotient case}
In this section, we study symplectic resolutions for quotient
symplectic varieties.  Let $V = \cit^{2n}$ and $G < Sp(2n)$ a finite
sub-group. For an element $g \in G$, we denote by $V^g$ the linear
subspace of points fixed by $g$. An element $g$ is called a {\em
symplectic reflection} if $codim\ V^g = 2$.
\begin{Thm}[Verbitsky \cite{Ver}]
Suppose that $V/G$ admits a symplectic resolution $\pi: Z \to
V/G$. Then $G$ is generated by symplectic reflections.
\end{Thm}
\begin{proof}
Let $G_0$ be the subgroup of $G$ generated by symplectic
reflections, then one has a natural map $\sigma: V/G_0 \to V/G$.
Let $F = \cup_{g \notin G_0} p_0(V^g)$, where $p_0: V \to V/G_0$
is the natural projection. Then $codim F \geq 4$ and $V/G_0 -F \to
V/G - \sigma(F)$ is a non-ramified covering of degree $\sharp
G/G_0$.

Let $Z_0 = V/G_0 \times_{V/G} Z$ be the fiber product and $\pi_0:
Z_0 \to V/G_0$ the projection to the first factor. Then $\sigma_0:
Z_0 - \pi_0^{-1}(F) \to Z - R$ is a non-ramified covering with
degree $\sharp G/G_0$, where $R = \pi^{-1}(\sigma(F))$.

Now the semi-smallness of $\pi$ implies $codim(R) \geq 2$. So for
the fundamental groups, one has $\pi_1(Z -R) \simeq \pi_1(Z) = {1}
$, where the last equality follows from the fact that any resolution
of $V/G$ is simply-connected. This shows that the non-ramified
covering $\sigma_0$ is of degree 1, thus $G = G_0$.
\end{proof}

\begin{Rque}
A classification of finite symplectic groups generated by
symplectic reflections is obtained in \cite{GS} (Theorem 10.1 and
Theorem 10.2) and also in \cite{Coh}.
\end{Rque}

Here are some examples of finite symplectic groups. One should
bear in mind that there are finite symplectic groups which are not
of this type. Let $L$ be a complex vector space and $G < GL(L)$ a
finite sub-group. Then $G$ acts on $L \oplus L^*$ preserving the
natural symplectic structure. Recall that an element $g \in G$ is
a {\em complex reflection} if $codim(L^g) = 1$. An element $g \in
G$ is a complex reflection if and only if it is a symplectic
reflection when considered as an element in $Sp(L\oplus L^*)$. Now
the precedent theorem implies:
\begin{Cor}[Kaledin \cite{Ka1}]
Suppose that $(L \oplus L^*)/G$ admits a projective symplectic
resolution $\pi: Z \to (L \oplus L^*)/G$, then $G$ is generated by
complex reflections.
\end{Cor}

The proof of D. Kaledin is different to the one presented above.
Here is an outline of his proof. One observes  that there exists a
natural $\cit^*$ action on $(L \oplus L^*)/G \simeq (T^* L)/G$. For
a symplectic resolution $\pi: Z \to (L \oplus L^*)/G$, one shows
that the $\cit^*$-action on $(T^* L)/G$ lifts to
 $Z$ in  such a  way that  $\pi$ is $\cit^*$-equivariant.
 In fact, this is a general fact for symplectic resolutions (see
 Lemma 5.12 \cite{GK}).

The key point is to show that if $\pi$ is furthermore projective,
then for every $x \in L/G \subset (T^*L)/G$, there exists at most
finitely many points in $\pi^{-1}(x)$ which are fixed by the
$\cit^*$-action on $Z$. The proof is based on the equation
$\lambda^* \Omega = \lambda \Omega$ and the semi-smallness of the
map $\pi$ (see Proposition 6.3 \cite{Ka1}).

Now since a generic point on $L/G$ is smooth, the map $\pi:
\pi^{-1}(L/G) \rightarrow L/G$ is generically one-to-one and
surjective, thus there exists a connected component $Y$ of fixed
points $Z^{\cit^*}$ such that $\pi: Y \rightarrow L/G$ is dominant
and generically one-to-one, which is also finite by the key point.
Now that $L/G$ is normal implies that $\pi: Y \rightarrow L/G$ is
in fact an isomorphism. Since $Z$ is smooth, $Z^{\cit^*}$ is a
union of smooth components, so $Y$ is smooth, thus $L/G$ is
smooth. Now a classical result then implies that $G$ is generated
by complex reflections.

This geometric approach can be developed further to obtain the
following theorem, which holds also for the more general case
$(T^*X)/G$ with $X$ a smooth variety and $G < Aut(X)$ a finite
sub-group.
\begin{Thm}[Fu \cite{Fu2}]
Let $L$ be a vector space and $G < GL(L)$ a finite sub-group.
Suppose we have a projective symplectic resolution $\pi: Z \to
(T^*L)/G$. Then:

(i).  $Z$ contains a Zariski open set $U$ which is isomorphic to
$T^*(L/G)$;

(ii). the restricted morphism $\pi: T^*(L/G) \to (T^*L)/G$ is the
natural one, which is independent of the resolution.
\end{Thm}

From what we have discussed above, there exists a connected
component $Y$ of $Z^{\cit^*}$ which is isomorphic to $L/G$. Notice
that the symplectic structure $\Omega$ on $Z$ satisfies $\lambda^*
\Omega = \lambda \Omega$ and $Y$ is of half dimension.  One deduces
that $Y$ is in fact a Lagrangian sub-manifold of $Z$. Now using a
classical result of A. Bialynicki-Birula (\cite{BB}), one can prove
that the attraction variety of $Y$ with respect to the
$\cit^*$-action is isomorphic to $T^*(L/G)$. Then one shows that
under this isomorphism, the morphism $T^*(L/G) \xrightarrow{\pi}
(T^*L)/G$ is the following natural one: take a point $[x] \in L/G$
and a co-vector $\alpha \in T_{[x]}^*(L/G)$. We define a co-vector
$\beta \in T_x^* L$ by $<\beta, v> = <\alpha, p_*(v) > $ for all $v
\in T_x L$, where $p: L \to L/G$ is the natural projection.  Then
$\pi ([x], \alpha) = [x, \beta]$.

A crucial question is how big the open set $U$ is in $Exc(\pi)$ and
how the rest part looks like. Some partial answer to this question
is obtained in \cite{Fu2}, where one needs the following version of
McKay correspondence:
\begin{Thm}[Kaledin \cite{Ka2}]
Let $V$ be a symplectic vector space and $G < Sp(V)$ a finite
sub-group.  Suppose we have a  symplectic resolution $\pi: Z \to
V/G$. Then there exists a basis (represented by maximal cycles of
$\pi$) of $H_{2k}(Z, \qit)$ indexed by the conjugacy classes of
elements $g \in G$ such that $codim \ V^g = 2k$.
\end{Thm}

\begin{Exam}
One special case is the following: let $\g$ be a complex semi-simple
Lie algebra and $\h$ a Cartan subalgebra. Let $G$ be the Weyl group
acting on $\h$. Then $W: = (\h \oplus \h^*)/G$ is a symplectic
variety. In the case of simple Lie algebras, it is proved in
\cite{GK} (Theorem 1.1) that $W$ admits a symplectic resolution if
and only if $\g$ is of type $A, B$ or $C$.

The case of type $A$ can be constructed as follows: let
$Hilb^n(\cit^2) \xrightarrow{\pi} (\cit^2)^{(n)}$ be the
Hilbert-Chow morphism and $\Sigma: (\cit^2)^{(n)} \to \cit^2$ the
sum map. Then $\pi_0: (\Sigma \circ \pi)^{-1}(0) \to \Sigma^{-1}(0)$
is a symplectic resolution of $\Sigma^{-1}(0)$, and $\Sigma^{-1}(0)$
is nothing but $(\h \oplus \h^*)/\mathcal{S}_n$, where $\h$ is a
Cartan subalgebra of $\mathfrak{sl}_{n}$.

\end{Exam}

However it is difficult to decide for which $G < Sp(V)$ the quotient
$V/G$ admits a symplectic resolution. The following problem is open
even when $dim(V) = 4$.
\begin{Pb}
(i). Classify finite sub-groups $G < Sp(2n, \cit)$ such that
$\cit^{2n}/G$ admits a symplectic resolution. (ii). Parameterizes
all symplectic resolutions of $V/G$.
\end{Pb}

Another obstruction to the existence of a symplectic resolution
makes use of the so-called Calogero-Moser deformation of $V/G$ ,
which is a canonical deformation of $V/G$ constructed by
Etingof-Ginzburg (see also \cite{GK}).
\begin{Thm}[Ginzburg-Kaledin \cite{GK}]
Suppose that $V/G$ admits a symplectic resolution, then a generic
fiber of the  Calogero-Moser deformation of $V/G$ is smooth.
\end{Thm}

In \cite{Gor}, I. Gordon proved that a generic fiber of the
Calogero-Moser deformation is singular for the quotient of some
symplectic reflection groups(Coxeter groups of type $D_{2n}, E, F,
H$ etc.). One may expect that his method can be used to obtain more
examples of $V/G$ which do not admit any symplectic resolution.

\section{Nilpotent orbit closures}
\subsection{symplectic resolutions}
As shown in Corollary \ref{nilpsym}, the normalization of a
nilpotent orbit closure is a symplectic variety. We will now discuss
their symplectic resolutions.
\begin{Prop}
Every nilpotent orbit closure in $\mathfrak{sl}(n+1)$ admits a
symplectic resolution.
\end{Prop}
\begin{proof}
Let $\0$ be the nilpotent orbit corresponding to the partition
$[d_1, \ldots, d_k]$. The dual partition is defined by $s_j = \sharp
\{i | d_i \geq j\}. $ The closure of $\0$  is
$$\overline{\0} = \{A \in \mathfrak{sl}(n+1) | \dim \ker A^j \geq \sum_{i=1}^j s_i \},
$$which is normal (\cite{KP}).

We define a  flag variety $F$ as follows:
$$F = \{(V_1, \ldots, V_l) | V_j \text{ vector space of dim} \sum_{i=1}^j s_i, V_j \subset V_{j+1} \text{ for all }
j\},$$ whose cotangent bundle $T^*F$ is isomorphic to the
coincidence variety $Z := \{(A, V_\bullet) \in \mathfrak{sl}(n+1)
\times F | A V_i \subseteq V_{i-1} \ \forall i\}$. The projection
from $Z$ to the first factor gives a morphism $\pi: T^*F \to
\overline{\0}$, which is in fact a resolution. Since $T^*F$ has
trivial canonical bundle, $\pi$ is a symplectic resolution of
$\overline{\0}$.
\end{proof}

One may wonder if every nilpotent orbit closure admits a symplectic
resolution. Unfortunately this is not the case, as shown by the
following:
\begin{Prop} \label{omin}
Let $\g$ be a simple Lie algebra. Then the closure
$\overline{\0}_{min}$ admits a symplectic resolution if and only if
$\g$ is of type $A$.
\end{Prop}
\begin{proof}
The Picard group of $\0_{min}$ is $\zit_2$ when $\g$ is of type $C$,
and is 0 if $\g$ is not of type $A$ or $C$. Thus
$\overline{\0}_{min}$ is in fact a normal $\qit$-factorial variety.
Now by the argument in the proof of Proposition \ref{purecodim2},
one sees that $\overline{\0}_{min}$ admits no symplectic resolution
if $\g$ is not of type $A$.
\end{proof}

Now the question is how to decide whether a nilpotent orbit
closure admits a symplectic resolution or not. If yes, can we find
 all of its symplectic resolutions? This question is answered by
the following
\begin{Thm}[Fu \cite{Fu1}]\label{nilp}
Let $\g$ be a semi-simple complex Lie algebra and $G$ its adjoint
group.  Let $\widetilde{\0}$ be the normalization of a nilpotent
orbit closure. Suppose that we have a symplectic resolution $\pi: Z
\to \widetilde{\0}$, then there exists a parabolic subgroup $P$ of
$G$ such that $Z$ is isomorphic to $T^*(G/P)$. Furthermore under
this isomorphism, the map $\pi$ becomes the moment map with respect
the action of $G$ (where $\g$ is identified with its dual via the
Killing form).
\end{Thm}

Recall that a parabolic subgroup $P$ is called a {\em polarization}
of $\0$ if $\overline{\0}$ is the image of the map $T^* (G/P) \to
\g$. Every parabolic subgroup is  a polarization of some nilpotent
orbit, but not every nilpotent orbit admits a polarization and those
admitting a polarization are called {\em Richardson orbits}.
\begin{Cor}
The normalization $\widetilde{\0}$ of a nilpotent orbit closure in a
semi-simple Lie algebra admits a symplectic resolution if and only
if (i). $\0$ is a Richardson nilpotent orbit; (ii). there exists a
polarization $P$ such that the moment map $T^*(G/P) \to
\widetilde{\0}$ is birational.
\end{Cor}

The key observation is that there exists a $\cit^*$ action on
nilpotent orbits, which follows directly from the Jacobson-Morozov
theorem (Theorem \ref{JM}). This $\cit^*$ (and $G$) action lifts not
only to the normalization $\widetilde{\0}$, but also to the
symplectic resolution $Z$. If we denote by $\Omega$ the symplectic
form on $Z$, then one feature of the $\cit^*$ action is $\lambda^*
\Omega = \lambda \Omega$. Together with the results of \cite{BB},
one shows that there exists an open set $U$ in $Z$ which is
isomorphic to $T^*Z_0$, where $Z_0$ is a connected component of
$Z^{\cit^*}$.

Now the action of $G$ on $Z$ restricts to an action on $Z_0$, which
is in fact transitive, thus $Z_0$ is isomorphic to $G/P$ for some
parabolic sub-group of $G$. Furthermore the restricted morphism of
$\pi$ to $U$ is in fact the moment map, which is a proper morphism,
thus $U$ is the whole of $Z$.

Using results in \cite{Hes} on polarizations of nilpotent orbits,
one obtains (see \cite{Fu1}) a classification of nilpotent orbit
closures of classical type whose normalization admits a symplectic
resolution .

This result can be generalized to odd degree coverings of
nilpotent orbits (see \cite{Fu3}), where an interesting phenomenon
appears: there exist some nilpotent orbits which admit some
symplectic resolutions, but not their coverings, and there exist
some nilpotent orbits which do not admit any symplectic
resolution, while some of their coverings do admit some symplectic
resolutions. A similar phenomenon appears also in \cite{KL}.

\subsection{birational geometry}

As shown in the precedent section, every symplectic resolution of a
nilpotent orbit closure  is of the form $T^*(G/P) \to
\widetilde{\0}$. However there can exist  several polarizations
which give different symplectic resolutions $T^*(G/P_i) \to
\widetilde{\0}, (i=1, 2).$ The birational geometry of the two
resolutions is encoded in the rational map $T^*(G/P_1) --\to
T^*(G/P_2)$. This section is to present several results on these
rational maps.

Consider the nilpotent orbit $\0= \0_{[2^k, 1^{n-2k}]}$ in
$\mathfrak{sl}_n$, where $2k < n$. The closure  $\overline{\0}$
admits exactly two symplectic resolutions given by $$T^*G(k,n)
\xrightarrow{\pi}  \overline{\0} \xleftarrow{\pi^+} T^*G(n-k, n),$$
where $G(k,n)$ (resp. $G(n-k, n)$) is the Grassmannian of $k$ (resp
$n-k$) dimensional subspaces
 in $\cit^n$. Let $\phi$ be the induced birational map
$T^*G(k,n) --\to T^*G(n-k,n)$. Then $\pi$ and $\pi^+$ are both
small and $\phi$ is a flop, which is called a {\em stratified
Mukai flop of type $A$.} These are the flops studied by E. Markman
in \cite{Mar}.

Let $\0$ be the orbit $\0_{[2^{k-1}, 1^2]}$ in $\mathfrak{so}_{2k}$,
where $k \geq 3$ is an odd integer. Let $G_{iso}^+(k), G_{iso}^-(k)$
be the two connected components of the orthogonal Grassmannian of
$k$-dimensional isotropic subspace in $\cit^{2k}$ (endowed with a
fixed non-degenerate symmetric form). Then we have two symplectic
resolutions $T^*G_{iso}^+(k) \to \overline{\0} \leftarrow
T^*G_{iso}^-(k)$.  This diagram is called a {\em stratified Mukai
flop of type $D$.}

Let $Z \xrightarrow{\pi} W \xleftarrow{\pi'} Z'$ be two resolution
of a variety $W$. Then the diagram is called a {\em locally trivial
family of stratified Mukai flops of type $A$ (resp. of type $D$)} if
there exists a partial open covering $\{U_i\}$ of $W$ which contains
the singular part of $W$ such that each diagram $\pi^{-1}(U_i) \to
U_i \leftarrow (\pi')^{-1}(U_i)$ is a trivial family of a stratified
Mukai flop of type $A$ (resp. of type $D$).

\begin{Thm}[Namikawa \cite{Nam4}]\label{nambir}
Let $\0$ be a nilpotent orbit in a classical complex Lie algebra and
$Z \to \widetilde{\0} \leftarrow Z'$ two symplectic resolutions.
Then the birational map $Z --\to Z'$ can be decomposed into finite
number of diagrams $Z_i \to W_i \leftarrow Z_{i+1} (i=1, \ldots,
k-1)$ with $Z_1  = Z$ and $Z_k = Z'$  such that each diagram is
locally a trivial family of stratified Mukai flops of type $A$ or of
type $D$.
\end{Thm}

 The proof in \cite{Nam4} consists of a case-by-case study, using
the classification of polarizations of a nilpotent orbit of
classical type in \cite{Hes}. The drawback is that nilpotent
orbits of exceptional type cannot be dealt with, since a
classification of polarization is not known in this case. However,
Y. Namikawa took another method in \cite{Nam5} by using Dynkin
diagrams instead of partitions to prove that a similar result
holds for nilpotent orbit closures in exceptional simple Lie
algebras, where when $\g$ is of type $E_6$, a new stratified Mukai
flop appears, and these are all flop types we need.

\section{Symplectic moduli spaces}

Consider a $K3$ or abelian surface $S$ endowed with  an ample
divisor $H$. Let $M_v$ be the moduli space of rank $r > 0$
$H$-semi-stable torsion free sheaves on $S$ with Chern class $(c_1,
c_2)$, where $v: = (r, c_1, c_2)$ is a Mukai vector. The open
sub-scheme $M_v^s$ in $M_v$ parameterizing $H$-stable sheaves is
smooth, whose tangent space at a point $[E]$ is canonically
isomorphic to $Ext^1_S(E, E)$.  The Yoneda coupling composed with
the trace map gives a bilinear form $Ext^1_S(E, E) \times Ext^1_S(E,
E) \to H^2(S, \0_S) = \cit$, which glues to  a symplectic form on
$M_v^s$ (\cite{Muk}, \cite{Bot}).

If $v$ is primitive, i. e. $g. c. d. (r, c_1 \cdot H, c_2) = 1$,
then $M_v^s = M_v$ is a smooth projective symplectic variety.
However for a multiple $v= m v_0$ of a primitive vector $v_0$ with
$m \geq 2$, the moduli space $M_v$ is singular. In the case of $m=
2$ and $\langle v_0, v_0\rangle = 2$, $M_v$ admits a unique
symplectic resolution constructed by O'Grady (\cite{O'G1},
\cite{O'G2}), where $\langle v_0, v_0\rangle = c_1^2 - 2 r c_2$ is
the Mukai pairing. What happens for other singular moduli spaces?

\begin{Thm}[Kaledin-Lehn-Sorger, \cite{KLS}]
Suppose that $H$ is $mv_0$-general, then the moduli space $M_{mv_0}$
is a projective symplectic variety which does not admit any
symplectic resolution if $m >2$ and $\langle v_0, v_0\rangle \geq 2$
or if $m \geq 2$ and $\langle v_0, v_0\rangle > 2$
\end{Thm}

In fact, they  proved that under the hypothesis of the theorem, the
moduli space $M_{mv_0}$ is locally factorial. Then the argument in
the proof of Proposition \ref{purecodim2} (see also Corollary 1.3
\cite{Fu1}) shows that $M_{mv_0}$ has no symplectic resolution,
since the codimension of the singular part is of codimension $\geq
4$.

For a $K3$ surface $S$, the case of $m =2$ has been proven in
\cite{KL}, and the case of $v= (2, 0, 2n) $ with $n \geq 3$ is
proved by Choy and Kiem in \cite{CK1}. For abelian surfaces, the
case of  $v= (2, 0, 2n)$ with $n \geq 2$ is proved in \cite{CK2}.
 The proof of Choy and Kiem
is based on another obstruction to the existence of a symplectic
resolution, which we present in the following.

The Hodge-Deligne polynomial of a variety $X$ is defined as  $$E(X;
u, v) = \sum_{p, q} \sum_{k \geq 0} (-1)^k h^{p,q}(H_c^k(X,\cit))
u^p v^q, $$ where $h^{p,q}(H_c^k(X,\cit))$ is the dimension of $(p,
q)$ Hodge-Deligne component in the $k$th cohomology group with
compact supports.

Let $W$ be a symplectic variety and $p: X \to W$ a resolution of
singularities such that the exceptional locus of $p$ is a divisor
whose irreducible components $D_1, \cdots, D_k$ are smooth with only
normal crossings. Then $K_X = \sum_i a_i D_i$ with $a_i \geq 0$,
since $W$ has only rational Gorenstein singularities.  For any
subset $J \subset I: = \{1, \cdots, k\}$, one defines $D_J = \cap_{j
\in J} D_j$, $D_\emptyset = X$ and $D_J^0 = D_J - \cup_{i \in I-J}
D_i$.Then the stringy E-function of $W$ is defined by:
$$E_{st}(W; u, v) = \sum_{J \subset I} E(D_J^0;u, v) \prod_{j \in J} \cfrac{uv-1}{(uv)^{a_j+1} -1}.$$

\begin{Thm}[Batyrev \cite{Bat}]
The stringy function is independent of the choice of a resolution.
For a symplectic resolution $Z \to W$, one has $E_{st}(W;u,v) =
E(Z;u,v)$, in particular, the stringy function is a polynomial.
\end{Thm}

In \cite{CK1}, they used Kirwan's resolution to calculate the
stringy function of $M_v$ and then proved that this function is not
a polynomial, thus $M_{(2,0, 2n)}$ admits no symplectic resolution
for $n \geq 3$. Similar method is used in \cite{CK2}.

One may wonder if this method can be used to prove the non-existence
of a symplectic resolution for some quotients $\cit^{2n}/G$.
Unfortunately, this does not work. In fact, it is shown in
\cite{Bat2} that $E_{st}(\cit^{2n} /G) = \sum_i C_i(G) (uv)^{2n-i}$,
where $C_i(G)$ is the number of conjugacy classes in $G$ whose fix
point is of codimension $2i$. For the minimal nilpotent orbit closure
$\overline{\0}_{min}$ in a simple complex Lie algebra of classical type,
one calculates that only for type $D$, the stringy Euler function of
$\overline{\0}_{min}$ is not a polynomial, thus $\overline{\0}_{min}$
admits no symplectic resolutions (compare Proposition \ref{omin}).

\section{Some conjectures}

This section is to list some unsolved conjectures on symplectic
resolutions.

\subsection{finiteness} \label{finite}
Let $W$ be a symplectic variety and $Z \xrightarrow{\pi} W
\xleftarrow{\pi^+} Z^+$ two resolutions. Then $\pi$ and $\pi^+$ are
said {\em isomorphic} if the rational map $\pi^{-1} \circ \pi^{+}:
Z^+ --\to Z$ is an isomorphism. $\pi$ and $\pi^+$ are said {\em
equivalent} if there exists an automorphism $\psi$ of $W$ such that
$\psi \circ \pi$ and $\pi^+$ are isomorphic.
\begin{Conj} (\cite{FN})
A symplectic variety has at most finitely many non-isomorphic
symplectic resolutions.
\end{Conj}

For nilpotent orbit closures, this conjecture is verified, since
there are only finitely many conjugacy classes of parabolic
subgroups in a semi-simple Lie group. It is proved in \cite{FN} that
a symplectic variety of {\em dimension 4} has at most finitely many
non-isomorphic {\em projective} symplectic resolutions. Some
quotient varieties are shown to admit at most one {\em projective}
symplectic resolution in \cite{FN}. Such examples include  symmetric
products of a smooth symplectic surface, the quotient $\cit^{2n}/G$,
where $G < Sp(2n)$ is a finite sub-group whose symplectic
reflections form a single conjugacy class. The next case to be
studied is $\cit^{2n}/G$ for a general $G$.

The proof of this conjecture can be divided into two parts: (i). a
symplectic variety can have at most finitely many non-equivalent
symplectic resolutions; (ii). almost all automorphism of a
symplectic variety can be lifted to any symplectic resolution.

\subsection{deformation}\label{deformation}

Recall that a deformation of a variety $X$ is a
 flat morphism
$\mathcal{X} \xrightarrow{p} S$ from a variety $\mathcal{X}$ to a
pointed smooth connected curve $0 \in S$ such that $p^{-1}(0) \cong
X$. Moreover, a deformation of a proper morphism $f: X \to Y$ is a
proper $S$-morphism $F: \mathcal{X} \to \mathcal{Y}$, where
$\mathcal{X} \to S$ is a deformation of $X$ and $\mathcal{Y} \to S$
is a deformation of $Y$.

Let $ X \xrightarrow{f} Y \xleftarrow{f^+} X^+$ be two proper
morphisms. One says that $f$ and $f^+$ are {\em deformation
equivalent} if there exists deformations of $f$ and $f^+$: $\X
\xrightarrow{F} \Y \xleftarrow{F^+} \X^+$ such that for any general
$s \in S$ the morphisms $\X_s \xrightarrow{F_s} \Y_s
\xleftarrow{F^+_s} \X^+_s$ are isomorphisms.

As to the relation between two symplectic resolutions, we have the
following:

\begin{Conj}(\cite{FN}, \cite{Ka3})\label{deform}
Suppose that we have two symplectic resolutions $\pi_i: Z_i \to W,
i=1,2, $
 then $\pi_1$ is deformation equivalent to $\pi_2$, and $Z_1$ is diffeomorphic to $Z_2$.
\end{Conj}

The motivation of this conjecture is the well-known theorem of D.
Huybrechts (\cite{Huy}), which says that two birational compact
hyperk\"ahler manifolds are deformation equivalent.  This
conjecture is true when $W$ is projective (\cite{Nam1}).  For
nilpotent orbit closures in a simple Lie algebra,  this conjecture
is shown to be true in  \cite{Nam4} (see \cite{Fu4} for a weaker
version). Under a rather restrictive additional assumption, this
conjecture is proved in \cite{Ka3}. We proved in \cite{Fu5} that
any two projective symplectic resolutions of $\cit^4/G$ are
deformation equivalent, where $G < Sp(4)$ is a finite subgroup.

\subsection{cohomology}\label{cohomology}
By Conjecture \ref{deform},   the cohomology ring $H^*(Z, \cit)$ of
a symplectic resolution $Z \to W$  is independent of the resolution,
in particular this invariant can be regarded as an invariant of $W$,
instead of $Z$. How can we recover this invariant from $W$?

In the case of quotient varieties $V/G$, there exists an orbifold
cohomology $H_{orb}^*(V/G, \cit)$ which is isomorphic as an algebra
to $H^*(Z, \cit)$ for a symplectic resolution $Z \to V/G$ (see
\cite{GK}). However for a general variety, the orbifold cohomology
is not defined.

There exists a natural cohomology on a symplectic variety $W$, the
Poisson cohomology $HP^*(W)$ (see \cite{GK}). Our hope is
\begin{Conj}
Let $Z \to W$ be a symplectic resolution, then $H^*(Z, \cit) \simeq
HP^*(W)$ as vector spaces. In particular, $HP^*(W)$ is
finite-dimensional.
\end{Conj}

This is true if $W$ is itself smooth. For a symplectic resolution $Z
\to V/G$, it is proved in \cite{GK} that $H^i(Z, \cit) \simeq
HP^i(V/G)$ for $i = 0, 1, 2$.

\subsection{derived equivalence}\label{derived}

As a special case of the Bondal-Orlov-Kawamata's (see \cite{Kaw})
conjecture, one has:

\begin{Conj} Suppose that we have two symplectic resolutions $Z_i \to W,
i=1,2, $ then there is an equivalence of derived categories of
coherent sheaves $D^b(Z_1) \sim D^b(Z_2)$.
\end{Conj}

This is shown to be true for four dimensional symplectic varieties
by Y. Kawamata and Y. Namikawa independently, using the work of J.
Wierzba and  J. Wisniewski (\cite{WW}). For a symplectic
resolution of a quotient variety $Z \to V/G$, it is shown in
\cite{BK} that there exists an equivalence of derived categories
$D^b(Coh(Z)) \simeq D^b(Coh(V)^G)$. In particular, the conjecture
is verified in this case.

The next case to be studied is nilpotent orbit closures in a
classical simple Lie algebra. By Theorem \ref{nambir}, this is
essentially reduced to prove the equivalence for the stratified
Mukai flops of type $A, D$ and $E$.

Recently, D. Kaledin proved this conjecture locally on $W$ in
\cite{Ka5}. Furthermore, he showed that if $W$ admits an expanding
$\cit^*$ action (such as nilpotent orbit closures), then the
conjecture is true. However it is not easy to compute the
equivalent functor in any explicit fashion, contrary to the case
done by Y. Kawamata and Y. Namikawa.

\subsection{birational geometry}\label{birational}
One way of constructing a symplectic resolution from another is to
perform Mukai's elementary transformations (\cite{Muk}), which can
be described as follows. Let $W$ be a symplectic variety and $\pi: Z
\to W$ a symplectic resolution. Assume that $W$ contains a smooth
closed subvariety $Y$ and $\pi^{-1}$ contains a subvariety $P$  such
that the restriction of $\pi$  to $P$ makes $P$ a $\pit^n$-bundle
over $Y$. If $codim(P) = n$, then we can blow up $Z$ along $P$ and
then blow down along the other direction, which gives another
(proper) symplectic resolution $\pi^+: Z^+ \to W$, provided that
$Z^+$ remains in our category of algebraic varieties. The diagram $Z
\to W \leftarrow Z^+$ is called  {\em  Mukai's elementary
transformation}
 (MET for short) over $W$ with center $Y$. A {\em MET
in codimension 2} is a diagram which becomes a MET after removing
subvarieties of codimension greater than 2. The following conjecture
is proposed in \cite{HY} (see also the survey \cite{Hu}).

\begin{Conj}(\cite{HY})
Let $W$ be a symplectic variety which admits two projective
symplectic resolutions $\pi: Z \to W$ and $\pi^+: Z^+ \to W$. Then
the birational map $\phi = (\pi^+)^{-1} \circ \pi: Z --\to Z^+$ is
related by a sequence of METs over $W$ in codimension 2.
\end{Conj}

Notice that since the two resolutions $\pi, \pi^+$ are both
crepant, the birational map $\phi$ is isomorphic in codimension 1.
This conjecture is true for four-dimensional symplectic varieties
by the work of Wierzba and  Wi\'sniewski (\cite{WW}) (while
partial results have been obtained in \cite{BHL}, see also
\cite{CMSB}). In \cite{Fu2}, this conjecture is proved for
nilpotent orbits in classical simple Lie algebras. For quotient
varieties $\cit^{2n}/G$, this conjecture is recently proved in
\cite{Fu5}. One hopes that this conjecture could serve to prove
the deformation equivalence.

\quad \\
C.N.R.S.,
Laboratoire J. Leray (Math\'ematiques)\\
 Facult\'e  des sciences \\
2, Rue de la Houssini\`ere,  BP 92208 \\
F-44322 Nantes Cedex 03 - France\\
\quad \\
fu@math.univ-nantes.fr

\end{document}